\documentclass[a4paper,12pt,legno]{article}
\usepackage[polish]{babel}
\usepackage[cp1250]{inputenc}
\usepackage[T1]{fontenc}
\usepackage{amsmath}
\usepackage{amsfonts}
\usepackage{amsthm}
\usepackage{mathrsfs}

\newtheorem{stw}[]{Proposition}
\newtheorem{lem}[]{Lemma}

\newtheorem{wn}{Corollary}
\newtheorem{defi}{Definition}

\newtheorem{uw}[]{Remark}

\newcommand{\Rset}{\mathbb{R}} 
\newcommand{\J}{\mathcal{J}}   
\newcommand{\I}{\mathcal{I}}   
\newcommand{\La}{\mathcal{L}}  

\newcommand{\nn}{\nonumber}

\newcommand{\be}{\beta}

\newcommand{\nbs}{\!\!\!\!\!}
\newcommand{\ul}{\limits^}
\newcommand{\dl}{\limits_}
\newcommand{\bl}{\bigl(}
\newcommand{\br}{\bigr)}
\newcommand{\Bl}{\Bigl(}
\newcommand{\Br}{\Bigr)}
\newcommand{\wh}{^{\widehat{}}}
\providecommand{\norm}[1]{\left\lVert#1\right\rVert}

\providecommand{\nor}[1]{|\!|#1|\!|}

\begin{document}

\title{Factorization property of generalized s-selfdecomposable measures and class  $L^f$ distributions$^1$}

\author{Agnieszka Czyzewska-Jankowska and Zbigniew J. Jurek\footnote{Corresponding Author. \newline $^1$Research funded by grant MEN Nr
1P03A04629, 2005-2008.}}

\date{Teor. Verojatn. Primen. 55(2010), 812-819; or Theory Probab. Appl. 55 (2011), 692-698.}

\maketitle
\begin{quote} \textbf{Abstract.} The method of \emph{random integral
representation}, that is, the method of representing  a given
probability measure as the probability distribution of some random
integral, was quite successful in the past few decades. In this note
we will find  such a representation for generalized
s-selfdecomposable  and selfdecomposable distributions that have the
\emph{factorization property}. These are the classes
$\mathcal{U}^f_{\be}$ and $L^f$, respectively



\emph{Abbrivated title:} Factorization property

\end{quote}

\medskip
\medskip
\medskip
In probability theory, from its very beginning, characteristic
functions (Fourier transforms) were used to describe measures and to
prove limiting distributions theorems. In the past few decades many
classes of probability measures (e.g. selfdecomposable measures ,
n-times selfdecomposable, s-selfdecomposable, type G distribution,
etc.) were characterized in terms of distributions of some random
integrals; cf. Jurek (1985, 1988) , Jurek and Vervaat (1983), Jurek
and Mason (1993), Jurek and Yor (2004), Iksanov, Jurek and Schreiber
(2004) and recently Aoyama and Maejima (2007). More precisely, for
each of those classes one integrates a fixed deterministic function
with respect to a class of L\'evy processes, with possibly a time
scale change.

Moreover, what we must emphasize here is that from the random
integral representations easily follow those in terms of
characteristic functions, and also one can infer from them new
convolution factorizations or decompositions. Thus the random
integral representations provide a new method in the area called
\emph{the arithmetic of probability measures}; cf. Cuppens (1975) or
Linnik and Ostrovskii (1977).

In this note we consider more specific situations. Namely, for a
convolution semigroup $\mathcal{C}$ of distributions of some random
integrals and  a measure $\mu\in\mathcal{C}$ we are interested in
decompositions of the form
\begin{equation}
\mu=\mu_1\ast\rho , \ \ \mu_1\in\mathcal{C},
\end{equation}
for some probability measure $\rho$ that is intimately related to
the measure $\mu_1$.

This paper was inspired by questions related to the class $L^f$ of
selfdecomposable measures having the so called \emph{factorization
property} that was introduced and investigated in Iksanov, Jurek and
Schreiber (2004).

Finally, let us note that the random integral representations for
classes $\mathcal{U}^f_{\be}$ (Corollary 1(a)) and $L^f$ (Corollary
3) provide more examples for the conjectured \ "meta-theorem"  \ in
\ \emph{The Conjecture} on \ www.math.uni.wroc.p/$\sim$zjjurek \, or
see Jurek (1985) and (1988).

\medskip
\medskip
\textbf{1. Notation and the results.} Our results are presented for
probability measures on Euclidean space $\Rset^d$. However, our
proofs are such that they hold true for measures on infinite
dimensional real separable Banach space $E$ with  the \emph{scalar
product}  replaced by the \emph{bilinear form} between $E^{\prime}
\times E$ and $\Rset$; $E^{\prime}$ denotes the topological dual of
$E$ and, of course, $(\Rset^d)^{\prime}=\Rset^d$; cf. Araujo-Gin\'e
(1980), Chapter III. In particular,  one needs to keep in mind
Remark 1, below.

\medskip
\medskip
Let $ID$ and $ID_{\log}$ denote all infinitely divisible probability
measures (on $\Rset^d$ or $E$) and those that integrate the
logarithmic function $\log(1+||x||)$, respectively. Let $Y_{\nu}(t),
t\ge0$ denote an $\Rset^d$ (or $E$) - valued L\'evy process, i.e., a
process with stationary independent increments, starting from zero,
and with paths that continuous from the right and with finite left
limits, such that $\nu$ is its probability distribution at time 1:
$\mathcal{L}(Y_{\nu}(1))=\nu$, where $\nu$ can be any $ID$
probability measure.

Throughout the paper $\mathcal{L}(X)$ will denote the probability
distribution of an $\Rset^d$-valued random vector (or a Banach space
E-valued random elements if the Reader is interested in that
generality).

 \begin{defi}
 For $\beta>0$ and a L\'evy process $Y_{\nu}$, let us define
\begin{equation}
\mathcal{J}^{\beta}(\nu):\,=\La\bl\int_0^1
t^{1/\beta}\;dY_{\nu}(t)\br=\La\bl\int_0^1
t\;dY_{\nu}(t^{\beta})\br, \ \
\mathcal{U}_{\beta}:\,=\mathcal{J}^{\beta}(ID).
\end{equation}
To the distributions from $\mathcal{U}_{\beta}$ we  refer to as
generalized s-selfdecomposable distributions.
\end{defi}
The classes $\mathcal{U}_{\beta}$ were already introduced in Jurek
(1988) as the limiting distributions in some schemes of summing
independent variables. The terminology has its origin in the fact
that distributions from the class $\mathcal{U}_{1} \equiv
\mathcal{U}$ were called \emph{s-selfdecomposable distribution} (the
"\emph{s}-" , stands here for \emph{the shrinking operations} that
were used originally in the definition of $\mathcal{U}$); cf. Jurek
(1985), (1988) and references therein.

\begin{stw} \textbf{A factorization of generalized s-selfdecomposable distribution.}
In order that a generalized s-selfdecomposable distribution $\mu
=\J^{\be}(\rho)$, from the class $\mathcal{U}_\be$,  convoluted with
its background measure $\rho$ is again in the class
$\mathcal{U}_\be$ it is sufficient and necessary that
$\rho\in\mathcal{U}_{2\be}.$

More explicitly,
\begin{equation}
[\,\J^\be(\rho)*\rho=\J^\be(\nu)\,] \Longleftrightarrow
[\,\rho=\J^{2\be}(\nu^{*\tfrac{1}{2}})\,]
\end{equation}

Furthermore, for each $\tilde \mu \in \mathcal{U}_\be$ there exists
a unique  $\tilde\rho\in \mathcal{U}_{2\be}$ such that
$\tilde\mu=\J^\be(\tilde\rho)* \tilde\rho$ and
$\J^{2\be}\bigl(\tilde\mu\bigr)=\J^{\be}\bigl((\tilde\rho)^{{\displaystyle
*}2}\bigr)$
\end{stw}
Let us denote by $\mathcal{U}_{\be}^f$ the class of generalized
s-selfdecomposable admitting \emph{the factorization property}, i.e,
 $\mu:=\J^\be(\rho)\in\mathcal{U}_{\be}$ has the factorization
property if $\J^\be(\rho)\ast\rho \in \mathcal{U}_{\be}$.
\begin{wn} For $\beta>0$ we have equalities
\begin{multline*}
(a) \ \
\mathcal{U}^f_{\be}=\mathcal{J}^{2\be}(\mathcal{U}_{\beta})=\mathcal{J}^{2\be}(\mathcal{J}^\be(ID))=
\\ = \{\La(\int^{1}_{0}(1-\sqrt{t})^{1/\be}\,dY_{\nu}(t)): \nu\in
ID\}. \qquad \qquad \qquad
\end{multline*}

(b) $\mathcal{U}_{\be}=\{\mathcal{J}^{\be}(\rho)\ast\rho : \rho\in
\mathcal{U}_{2\be}\}$.
\end{wn}
Taking in Proposition 1 $\be=1$ we get the  following
\begin{wn}{\textbf{Factorization of s-selfdecomposable distributions. }}
An s-selfdecomposable distribution $\mu =\J(\rho)$ convoluted with
$\rho$ is again s-selfdecomposbale if and only if
$\rho\in\mathcal{U}_2$. Thus we have
$\mathcal{U}^f=\mathcal{J}^{2}(\mathcal{U})$.

More explicitly
\begin{equation}
[\,\J(\rho)*\rho=\J(\nu)\,] \Longleftrightarrow [\,\rho=
\J^{2}\bigl(\nu^{{\displaystyle *}\tfrac{1}{2}}\bigr)\,].
\end{equation}
Moreover, for each $\tilde\mu \in \mathcal{U}$ there exist a unique
$\rho\in \mathcal{U}_2$ such that
$\tilde\mu=\J(\tilde\rho)*\tilde\rho$ and
$\J^{2}\bigl(\tilde\mu\bigr)=\J\bigl((\tilde\rho)^{*2})\bigr)$.
Consequently,  $\mathcal{U}=\{\mathcal{J}^{2}(\rho)\ast\rho :
\rho\in \mathcal{U}\}$.
\end{wn}

Following Jurek-Vervaat (1983) or Jurek (1985) we recall the
following

\begin{defi} For a measure $\nu\in ID_{\log}$ and a L\'evy process
$Y_{\nu}$ let us define
\begin{equation}
 \mathcal{I}(\nu):= \mathcal{L}\bigl(\int_0^{\infty}e^{-s}\,
d\,Y_{\nu}(s)\bigr), \ \  \  L:=\mathcal{I}(ID_{\log})
\end{equation}
and distributions from $L$ are called selfdecomposable or L\'evy
class L distributions.
\end{defi}

In classical probability  theory the selfdecomposability ( or in
other words, the L\'evy class $L$ distributions) is usually defined
via some decomposability property or by scheme of limiting
distributions. However,  since Jurek-Vervaat (1983) we know that the
class $L$ coincides with the class of distributions of random
integrals given in (5) and  thus it is used in this note as its
definition.

\medskip
Before going further, let us recall the following example that led
to, and justified interest in, that kind of
investigations/factorizations.

\medskip
\textbf{Example.} For two dimensional Brownian motion  $\textbf{B}_t
:= (B^1_t, B^2_t)$, the process
\[
\mathcal{A}_t:=\int_0^tB^1_s\,dB^2_s-B^2_s\,dB^1_s,\ \ \  t>0,
\]
called \emph{L\'evy's stochastic area integral}, admits the
following factorization
\begin{equation}
\chi(t):=E[ e^{it\mathcal{A}_u}|\textbf{B}_u=(\sqrt{u},
\sqrt{u})]=\frac{tu}{\sinh tu}\cdot \exp[-(tu\,\cosh tu -1)],
\end{equation}
cf. P. L\'evy (1951) or Yor (1992), p. 19.

Iksanov-Jurek-Schreiber (2004), p. 1367, proved that the
factorization (6) may be interpreted as follows: if $\nu$ is the
probability measure with the characteristic function $t\to
\exp[-(tu\,\cosh tu-1)]$ then $\mathcal{I}(\nu)$  has the
characteristic function $t\to \frac{tu}{\sinh tu}$, and also
\begin{equation}
\mathcal{I}(\nu)\ast \nu = \mathcal{I}(\rho), \ \ \mbox{ for some} \
\ \rho\in ID_{\log};
\end{equation}
i.e., $\mathcal{I}(\nu)$  is selfdecomposable and when convoluted
with its background driving probability measure $\nu$ we again  get
a distribution from the class $L$.

Let us note that the convolution factorizations (7), (3) and (4) are
of the form described in (1), with different semigroups
$\mathcal{C}$.

\medskip
\begin{stw}{\textbf{Random integral representation
of
$\boldsymbol{\I\bigl(\J^\be\bigl(\textsl{ID}_{\log}\bigr)\bigr)}$}}.
\\ For $\nu\in\textsl{ID}_{\log}$ and $\be>0$

\begin{equation}
\I\bigl(\J^\be\bigl(\nu\bigr)\bigr)
=\La\bigl(\int_{0}^{\infty}{\displaystyle e^{\displaystyle
-s}}\,dY_{\nu}\bl \sigma_{\be}(s)\br\bigr) ,
\end{equation}
where $Y_{\nu}(t), t\ge 0$ is a L\'evy process such that
$\mathcal{L}(Y_{\nu}(1))=\nu$ and the deterministic inner clock
$\sigma_{\beta}$ is given by $
 \sigma_{\be}(s):= s+\tfrac{1}{\be}e^{-\be s}-\tfrac{1}{\be}, \
 s\ge0$.
 \end{stw}

From Proposition 1 (ii) in Iksanov-Jurek-Schreiber (2004) and taking
$\be=1$ in Proposition 2 we get

\begin{wn} For the class, $L^f$, of selfdecomposable distributions with factorization property,  we  have the
 following random integral representation
\begin{equation}
 L^f = \big\{\La\bigl(\int_{0}^{\infty} e^{-s}\,dY_{\nu}(s +e^{-s}-1))\,: \ \nu\in
 ID_{\log}\big\}.
\end{equation}
\end{wn}

\medskip
\textbf{2. Proofs.} For a probability Borel measures $\mu$ on
$\Rset^d$, its \emph{characteristic function} $\hat{\mu}$ is defined
as
\[
\hat{\mu}(y):=\int_{\Rset^d} e^{i<y,x>}\mu(dx), \ y\in\Rset^d,
\]
where $<\cdot,\cdot>$ denotes the scalar product;  (in case one
wants to have results on Banach spaces $<\cdot,\cdot>$ is the
bilinear form on $E^{\prime}\times E$ and $y\in E^{\prime}$).

\noindent Recall that for infinitely divisible measures $\mu$ their
characteristic functions admit the following L\'evy-Khintchine
formula
\begin{multline}
\hat{\mu}(y)= e^{\Phi(y)}, \ y \in \Rset^d, \ \ \mbox{and the
exponents $\Phi$ are of the form} \\  \Phi(y)=i<y,a>-
\frac{1}{2}<y,Sy>  +  \qquad \qquad \qquad \\ \int_{\Rset^d
\backslash \{0 \}}[e^{i<y,x>}-1-i<y,x>1_B(x)]M(dx),
\end{multline}
where $a$ is  a \emph{shift vector}, $S$ is a \emph{covariance
operator} corresponding to the Gaussian part of $\mu$ and $M$ is a
\emph{L\'evy spectral measure}. Since there is a one-to-one
correspondence between a measure $\mu \in ID$ and the triples $a$,
$S$ and $M$ in its L\'evy-Khintchine formula (10) we will write
$\mu=[a,S,M]$. Finally, let recall that
\begin{equation}
M \ \mbox{is L\'evy spetral measure on $\Rset^d$ iff} \ \
\int_{\Rset^d}\min(1, ||x||^2)M(dx)<\infty
\end{equation}
(For infinite divisibility of probability measures on Banach spaces
we refer to the monograph by Araujo-Gin\'e (1980), Chapter 3,
Section 6, p. 136. Let us stress that the characterization (11), of
L\'evy spectral measures, is in general NOT true in infinite
dimensional Banach spaces ! However, it holds true in Hilbert
spaces; cf. Parthasarathy (1967), Chapter VI, Theorem 4.10.)

\medskip
Before proving Proposition 1, let us note the following auxiliary
facts.
\begin{lem}
(a) For the mapping $\J^\be$  and $\nu\in ID$  we have
\begin{equation}
 \widehat{\mathcal{J}^{\beta}(\nu)}(y)=\exp\int_0^1\log
\widehat{\nu}(t^{1/\beta}y)\,dt =
\exp\mathbf{E}[\log\widehat{\nu}(U^{1/\beta}y)], \ \ y\in\Rset ^d \
(\mbox{or}\ E^{\prime}).
\end{equation}
and $U$ is a random variable uniformly distributed over the unit
interval $(0,1)$.

(b) The mapping $\mathcal{J}^{\beta}$ is one-to-one. More explicitly
we have that
\begin{equation}
\frac{d}{ds}[s\log
\widehat{\mathcal{J}^{\beta}(\nu)}(s^{1/{\beta}}y)]|_{s=1}=
\log\hat{\nu}(y), \ \   \mbox{for all}  \ y\in \Rset^d \ (\mbox{or}\
E^{\prime}).
\end{equation}

(c) The mappings $\mathcal{J}^\beta, \beta >0$ commute, i.e., for
$\beta_1, \beta_2 >0$ and $\nu \in ID$, \ \ \ \ \ \ \ \ \ \
$\mathcal{J}^{\beta_1}(\mathcal{J}^{\beta_2}(\nu))=\mathcal{J}^{\beta_2}(\mathcal{J}^{\beta_1}(\nu))$.

(d) For probability measures $\nu_1, \nu_2$ and $c>0$ we have that
\begin{equation}
\mathcal{J}^{\beta}(\nu_1\ast\nu_2)= \mathcal{J}^{\beta}(\nu_1)
\ast\mathcal{J}^{\beta}(\nu_2); \ \ (\mathcal{J}^{\beta}(\nu))^{\ast
c}=\mathcal{J}^{\beta}(\nu ^{\ast c})
\end{equation}

(e) For $\beta >0$ and $\rho\in ID$ we have the identity
\begin{equation}
\mathcal{J}^{2\beta}(\mathcal{J}^{\beta}(\rho) \ast \rho) =
\mathcal{J}^{\beta}(\rho^{\ast 2})
\end{equation}
\end{lem}
\emph{Proof of Lemma 1.} Part (a) follows from the definition of the
random integrals and is a particular form (take matrix $Q=I$) of
Theorem 1.3 (a) in Jurek (1988).

\noindent For the claim (b) note that for each fixed $y$ we have
\[
\log\widehat{\mathcal{J}^{\beta}(\nu)}(s^{1/{\beta}}y)=s^{-1}\int_0^s\log\hat{\nu}(r^{1/{\beta}}y)dr,
\ \ s \in \Rset^+.
\]
This gives the  formula in (b), similarly as in Jurek (1988), p.
484. Equalities in (c) and (d) are also consequences of (a); cf.
Jurek(1988), Theorem 1.3 (a) and (c).

Finally, for the identity in (e) note, using (14) that
\begin{multline}
\log\Big(\mathcal{J}^{2\beta}\big(\mathcal{J}^{\beta}(\rho)\ast\rho\big)\Big)^{\widehat{}}(y)=
\int_0^1\log
\big(\mathcal{J}^{\beta}(\rho)\ast\rho\big)\Big)^{\widehat{}}(s^{1/2\beta}y)
=
\\ \int_0^1 \int_0^1\log\hat{\rho}(t^{1/\beta}s^{1/2\beta}y)dt\,ds +
\int_0^1\log\hat{\rho}(s^{1/2\beta}y)ds  \ \ \ \ (\mbox{put} \ t^2s=:u) \ \ \ \ \ \\
= \int_0^1 1/2 \int_0^s\log\hat{\rho}(u^{1/2\beta}y)
(us)^{-1/2}du\,ds + \int_0^1\log\hat{\rho}(s^{1/2\beta}y)ds \\=
\int_0^1\log\hat{\rho}(u^{1/2\beta}y)\,u^{-1/2} \big( 1/2\int_u^1
s^{-1/2}ds \big)\,du + \int_0^1\log\hat{\rho}(s^{1/2\beta}y)ds= \\
\int_0^1 u^{-1/2}\log\hat{\rho}(u^{1/2\beta}y)du= 2 \int_0^1
\log\hat{\rho}(u^{1/2\beta}y)d(u^{1/2})= \\ \int_0^1\log
\hat{\rho^{\ast 2}}(s^{1/\beta}y)ds =  \log\,
(\mathcal{J}^\beta(\rho^{\ast 2})){\hat{}}\,(y), \ \ \ \ \ \ \
\end{multline}
which completes the proof of Lemma 1.

\medskip
\emph{Proof of Proposition 1.} Suppose we have that
$\J^\be(\rho)*\rho=\J^\be(\nu)$.  Then by the properties described
in Lemma 1,
\[
\J^{\be}\bigl(\J^{2\be}(\nu)\bigr)=\J^{2\be}\bigl(\J^\be(\nu)\bigr)=\J^{2\be}\bigl(\J^\be(\rho)*\rho\bigr)=\J^\beta
(\rho^{\ast 2}),
\]
and hence $\rho^{\ast 2}= \J^{2\be}(\nu)$, i.e.,
$\rho=(\J^{2\be}(\nu))^{\ast 1/2}= \J^{2\be}(\nu^{\ast1/2})$, which
proves the necessity. The converse claim also follows from the above
reasoning.

For the last part, let us  note that  if $\tilde\mu=\J^{\beta}(\nu)
\in \mathcal{U}_{\beta}$ then taking
$\rho:=\J^{2\beta}(\nu^{\ast1/2})\in \mathcal{U}_{2\be}$ one gets
the required equality.

\medskip
\emph{Proof of Corollary 1.} Note that $\nu=\mathcal{J}^\be\in
\mathcal{U}^f_\be \ \mbox{iff} \ \J^\be(\rho)\ast\rho \in
\mathcal{U}_\be \ \ \mbox{iff} \ \rho \in \mathcal{U}_{2\be}$, by
(3) in Proposition 1. Last equality is from the Example \textbf{(a)}
from Czy�ewska-Jankowska and Jurek (2008). Similarly one gets part
(b) using Proposition 1 and Lemma 1 (e).

\medskip
Proposition 1 can be expressed in terms of characteristic
functions as follows:
\begin{wn}
In order that
\[
\exp{\int^{1}_{0}\log{\hat{\rho}\left(
t^{1/\be}y\right)}dt}\;\cdot\;\hat{\rho}\left(y\right)=
\exp{\int^{1}_{0}\log{\hat{\nu}\left( t^{1/\be}y\right)}dt},\qquad
y\in \Rset^d \ (\mbox{or}\ E')
\] for some $\mu$ and $\rho$ in ID it is necessary and
sufficient that
\[\hat{\rho}\left(y\right)=\exp{\int^{1}_{0}\tfrac{1}{2}\log{\hat{\nu}\left( t^{1/(2\be)}y\right)}dt};\]
\end{wn}

or in terms of the L\'evy spectral measures as:

\begin{wn}
In order to have the equality
\[\int^{1}_{0}M(t^{-1/\be}A)\;dt+M(A)=\int^{1}_{0}G(t^{-1/\be}A)\,dt,\qquad \textrm{for each Borel }A\in\mathcal{B}_0 ,\]
for some L\'evy spectral measures $M$ and $G$, it is necessary and
sufficient  that
\[M(A)=\int^{1}_{0}\tfrac{1}{2}G(t^{-1/(2\be)}A)\,dt, \qquad
\textrm{for each} \ \ A\in\mathcal{B}_0 ,
\]
\end{wn}
because if $\rho=[a, S, M]$ then the left hand side in the Corollary
is the L\'evy spectral measure of $\J^\be(\rho)*\rho$.

\medskip
For references let state the following
\begin{lem}
(i) If $\nu = [a,R,M]$ and
$\J^\be(\nu)=[a^{(\be)},R^{(\be)},M^{(\be)}]$ then
\begin{multline*}
a^{(\be)}:=\tfrac{\be}{(1+\be)}\;a+\int^{1}_{0}t^{1/\be}
\!\!\int_{\left\{1<\nor{x}\leq t^{-1/\be}\right\}}x\;M(dx)\;dt \\
=\frac{\be}{\be+1}\,(a+\int_{(||x||>1)}x\,||x||^{-1-\be}M(dx)\,); \
\ \ \ \ \
R^{(\be)}:=\tfrac{\be}{2+\be}\,R;  \qquad \qquad   \\
M^{(\be)}(A):=\int^{1}_{0}T_{t^{1/\be}}\;M(A)\;dt,\ \ \textrm{for
each}\,\,A\in\mathcal{B}_0 . \qquad \qquad \qquad \qquad \qquad
\qquad \qquad
\end{multline*}
(ii) For $\be>0$, we have that  $\J^\be(\nu)\in ID_{\log}$ if and
only if $\nu\in ID_{\log}$.
\end{lem}
\emph{Proof of Lemma 2.} (i) \ Uniqueness of the triplets: \emph{a
shift vector a , Gaussian covariance R and L\'evy spectral measure
M} in the L\'evy-Khintchine formula and equation (12) in Lemma 1
give the expressions for $a^{(\be)}, R^{(\be)}$ and for $M^{(\be)}$;
for details cf. formulas  (1.10), (1.11) and (1.12) in Jurek (1988),
with the matrix $Q=I$.

For part (ii), note that since we have
\begin{gather*}
\int\limits_{\left\{\norm{x}>1\right\}}\!\!\!\!\!\log{\nor{x}}\,M^{(\be)}(dx)=
\int\limits^{1}_{0}\!\!\int\limits_{\left\{\nor{x}>1\right\}}\nbs\log{\nor{x}}\,T_{t^{1/\be}}M(dx)\;dt=\nn\\
=\int\ul{1}\dl{0}\!\!\int\dl{\left\{\nor{t^{1/\be}x}>1\right\}}\nbs\nbs\log{\nor{t^{1/\be}x}}\,M(dx)\,dt
=\int\ul{1}\dl{0}\!\!\int\dl{\left\{\nor{x}>\tfrac{1}{t^{1/\be}}\right\}}\nbs\nbs\log{(t^{1/\be}\nor{x})}\,M(dx)\,dt=
\end{gather*}
\begin{gather*}
=\!\!\!\int\dl{\left\{\nor{x}>1\right\}}\int\ul{1}\dl{\nor{x}^{-1/\be}}\nbs\log{(t^{1/\be}\nor{x})}\,dt\,M(dx)
=\!\!\!\int\dl{\left\{\nor{x}>1\right\}}\nbs\tfrac{1}{\nor{x}^\be}\!\!\int\ul{\nor{x}}\dl{\nor{x}^{1-1/\be^2}}\nbs\!\!
\be w^{\be-1}\,\log{w}\,dw\,M(dx)=\nn\\
=\!\!\!\int\dl{\left\{\nor{x}>1\right\}} \nbs\tfrac{1}{\nor{x}^\be}
\biggl[w^\be\log{w}-\tfrac{1}{\be}w^\be\biggr|^{w=\nor{x}}_{w=\nor{x}^{1-1/\be^2}}\biggr]\,M(dx)=\nn\\
=\int\dl{\left\{\nor{x}>1\right\}}\nbs\log{\nor{x}}\,M(dx)-
\nbs\int\dl{\left\{\nor{x}>1\right\}}[ \tfrac{1}{\be}
+\tfrac{1}{\nor{x}^{1/\be}}\bl(1-\tfrac{1}{\be^2})\log{\nor{x}}-\tfrac{1}{\be}\br]\,M(dx)\;\nn
\end{gather*}
and the last integral is finite (the integrand function is bounded
on $(||x||>1)$ and L\'evy spectral measures M are finite on the
complements of all neighborhoods of zero; comp. (11)), therefore
from the above we conclude that
\[
[\,\int\limits_{\left\{\norm{x}>1\right\}}\!\!\!\!\!\log{\nor{x}}\,M^{(\be)}(dx)
< \infty] \ \ \mbox{iff} \ \ [
\int\limits_{\left\{\norm{x}>1\right\}}\!\!\!\!\!\log{\nor{x}}\,M(dx)
< \infty].
\]
But since the function $u\to \log(1+u)$, for $u>0$, is sub-additive
therefore we may apply Proposition 1.8.13 in Jurek-Mason (1993) and
infer the claim (ii). This completes the proof of Lemma 2.

\medskip
\emph{Proof of Proposition 2.}  If  $\nu\in ID_{\log}$  then, by
Lemma 2, $\mathcal{J}^{\beta}(\mu)\in ID_{\log}$ and thus the
improper random integral \
$\int_0^{\infty}e^{-s}dY_{\mathcal{J}^\be(\nu)}(s)$ \ converges (is
well-defined) almost surely (in probability and in distribution);
cf. Jurek-Vervaat (1983), Lemma 1.1 or  Jurek (1985). Hence and
Lemma 1(a) we get that
\begin{multline}
\log{\bl\I\left(\J^\be\left(\nu\right)\right)}\br\wh\left(y\right) =
\int_0^{\infty}\log \widehat{\mathcal{J}^{\be}(\nu)}(e^{-s}y)ds =
\int_0^{\infty}\int_0^1 \log\hat{\nu}(v^{1/\be}e^{-s}y) dv ds \\
= \int_0^1\int_0^{v^{1/\be}}\log \hat{\nu}(uy)u^{-1}du= \int_0^1
(\int_{u^\be}^1dv)\log \hat{\nu}(uy)u^{-1}du dv =\\
\int^{1}_{0}\log{\hat{\nu }\left(
uy\right)}\bigl(u^{-1}-u^{\be-1}\bigr)\,du
=\int^{\infty}_{0}\log{\hat{\nu}\bigl(e^{\displaystyle-s}y\bigr)}\bigl(1-e^{-\be s}\bigr)\,ds=\nn\\
= \int^{\infty}_{0}\log{\hat{\nu }\bigl(e^{{\displaystyle
-s}}y\bigr)\,d\sigma_{\be}(s)}.
\end{multline}

On the other hand, the random integral
\[
\int_{0}^{\infty} e^{-s}\,dY_{\nu}(\sigma_{\be}(s)):\,=\,\lim_{b\to
\infty} \int_0^b e^{-s}\,dY_{\nu}(\sigma_{\be}(s))\ \mbox{exists in
distribution},
\]
(or in probability or almost surely) because the function
\begin{multline*}
y\to \lim_{b\to \infty} {\Bl{\La\bigl(\int_{0}^b{e^{-s}}\,dY_{\nu}(\sigma_{\be}(s))}\Br\wh(y)}\\
= \lim_{b\to
\infty}\exp\int^{b}_{0}\log{\hat{\nu}(e^{-s}}y)\,d\sigma_{\be}(s)=
\exp\int^{\infty}_{0}\log{\hat{\nu}(e^{-s}}y)\,d\sigma_{\be}(s),
\end{multline*}
is a characteristic function. Moreover, we have that
\[
\I\bigl(\J^\be\bigl(\nu\bigr)\bigr)=\La\bigl(\int_{0}^{\infty}{
e^{-s}}\,dY_{\nu}\bl \sigma_{\be}(s)\br\bigr),
\]
which completes a proof of Proposition 2.

\medskip
\begin{uw}
Our argument above is valid for infinite dimensional Banach spaces,
although one should be aware that in that generality convergence of
characteristic functions to a characteristic function does not
guarantee weak convergence of corresponding distributions (
probability measures); cf. Araujo-Gine (1980), Theorem 4.19 on p.
29.
\end{uw}
\emph{Proof of Corollary 3.} Recall that by definition
$L^f=\{\mathcal{I}(\mu): \mathcal{I}(\mu)\ast\mu\in L\}$. However,
in view of Proposition 1 (ii) in Iksanov-Jurek-Schreiber (2004) we
have $L^f=\mathcal{I}(\mathcal{J}(ID_{\log})$. Consequently, taking
$\be=1$ in Proposition 2 we get the corollary.

\medskip
\medskip
\medskip
\medskip
\begin{center}
\textbf{References}
\end{center}

\medskip
\noindent [1] T. Aoyama and M. Maejima (2007). Characterizations of
subclasses otf type G distributions on $\Rset^d$ by stochastic
random integral representation,\emph{Bernoulli}, vol. 13, pp.
148-160.

\noindent [2] A. Araujo and E. Gine (1980). \emph{The central limit
theorem for real and Banach valued random variables.} John Wiley \&
Sons, New York.

\noindent [3] R. Cuppens (1975). \emph{Decomposition of multivariate
probabilities.} Academic Press, New York.

\noindent [4] A. Czyzewska-Jankowska and Z. J Jurek (2008). A note
on a composition of two random integral mappings $\J^\be$ and some
examples, preprint.

\noindent [5] A. M. Iksanov, Z. J. Jurek and B. M. Schreiber (2004).
A new factorization property of the selfdecomposable probability
measures, \emph{Ann. Probab}. vol.~32, Nr 2, str. 1356-1369.

\noindent [6] Z. J. Jurek (1985). Relations between the
s-selfdecomposable and selfdecomposable measures. \emph{Ann.
Probab.} vol.13, Nr 2, str. 592-608.

\noindent [7] Z. J. Jurek  (1988). Random Integral representation
for Classes of Limit Distributions Similar to Lavy Class $L_{0}$,
 \emph{Probab. Th. Fields.} 78, str. 473-490.

\noindent [8] Z. J. Jurek and J. D. Mason (1993).
\emph{Operator-limit distributions in probability theory.} John
Wiley \&Sons, New York.

\noindent [9] Z. J. Jurek and W. Vervaat (1983). An integral
representation for selfdecomposable Banach space valued random
variables, \emph{Z. Wahrscheinlichkeitstheorie verw. Gebiete}, 62,
pp. 247-262.

\noindent [10] Z. J. Jurek and M. Yor (2004). Selfdecomposable laws
associated with hyperbolic functions, \emph{Probab. Math. Stat.} 24,
no.1, pp. 180-190.

\noindent[11] P. L\'evy (1951). Wiener's random functions, and other
Laplacian random functions, \emph{Proc. Second Berkeley Symposium
Math. Statist. Probab.} str.~171-178. Univ. California Press,
Berkeley.

\noindent[12] Ju. V. Linnik and I. V. Ostrovskii (1977).
\emph{Decomposition of Random Variables and Vectors.} American
Mathematical Society, Providence, Rhode Island.

\noindent[13] K. R. Parthasarathy (1967). \emph{Probability measures
on metric spaces}. Academic Press, New York and London.

\medskip
\noindent
Institute of Mathematics \\
University of Wroc\l aw \\
Pl.Grunwaldzki 2/4 \\
50-384 Wroc\l aw, Poland \\
e-mail: zjjurek@math.uni.wroc.pl \ \ or \ \ czyzew@math.uni.wroc.pl \\
www.math.uni.wroc.pl/$^{\sim}$zjjurek

\end{document}